# CONSISTENCY OF BAYESIAN PROCEDURES FOR VARIABLE SELECTION

By George Casella,[1] F. Javier Girón,
M. Lina Martínez and Elías Moreno[2]

*University of Florida, University of Málaga, University of Málaga and University of Granada*

It has long been known that for the comparison of pairwise nested models, a decision based on the Bayes factor produces a consistent model selector (in the frequentist sense). Here we go beyond the usual consistency for nested pairwise models, and show that for a wide class of prior distributions, including intrinsic priors, the corresponding Bayesian procedure for variable selection in normal regression is consistent in the entire class of normal linear models. We find that the asymptotics of the Bayes factors for intrinsic priors are equivalent to those of the Schwarz (BIC) criterion. Also, recall that the Jeffreys–Lindley paradox refers to the well-known fact that a point null hypothesis on the normal mean parameter is always accepted when the variance of the conjugate prior goes to infinity. This implies that some limiting forms of proper prior distributions are not necessarily suitable for testing problems. Intrinsic priors are limits of proper prior distributions, and for finite sample sizes they have been proved to behave extremely well for variable selection in regression; a consequence of our results is that for intrinsic priors Lindley's paradox does not arise.

**1. Introduction.** Bayesian estimation of the parameters of a given sampling model is, under wide conditions, consistent. That is, the posterior probability of the parameter is concentrated around the true value as the sample size increases, assuming that the true value belongs to the parameter space being considered. The case where the dimension of the parameter

Received November 2006; revised December 2007.

[1]Supported by NSF Grants DMS-04-05543, DMS-06-31632 and SES-0631588.

[2]Supported by Ministerio de Ciencia y Tecnología, Grant SEJ-65200 and Junta de Andaluciá Grant SEJ-02814.

space is infinite can be an exception [see Diaconis and Friedman (1986) for examples of inconsistency of Bayesian methods].

When several competing models are deemed possible, so that we have uncertainty among them, consistency of a Bayesian model selection procedure is much more involved. For instance, it is well known that improper priors for the model parameters cannot be used for computing posterior model probabilities. Therefore, the priors need be either proper or limits of sequences of proper priors. Furthermore, not every limit of proper priors is appropriate for a Bayesian model selection.

The so-called Lindley paradox is an example of this [Lindley (1957) and Jeffreys (1967)]; it shows that when testing a point null hypothesis on the normal mean parameter we always accept the null if a conjugate prior is considered on the alternative and the variance of this conjugate prior goes to infinity. As Robert (1993) has pointed out, this is not a mathematical paradox since the prior sequence is giving less and less mass to any neighborhood of the null point as the prior variance goes to infinity. However, an important consequence of the paradox is that some limiting forms of proper priors might not be suitable for testing problems as they could provide inconsistency of the corresponding Bayes factors. We remark that intrinsic priors are limits of sequences of proper priors [Moreno, Bertolino and Racugno (1998)] and for finite sample sizes an intrinsic Bayesian analysis has been proved to behave extremely well for variable selection in regression [Casella and Moreno (2006), Girón, Moreno and Martínez (2006) and Moreno and Girón (2008)]. Consequently, showing that the Lindley paradox does not occur when using intrinsic priors is an important point.

For nested models and proper priors for the model parameters, the consistency of the Bayesian pairwise model comparison is a well established result [see O'Hagan and Forster (2004) and references therein]. Assuming that we are sampling from one of the models, say $M_1$, which is nested in $M_2$, consistency is understood in the sense that the posterior probability of the true model tends to 1 as the sample size tends to infinity. We observe that the posterior probability is defined on the space of models $\{M_1, M_2\}$. An equivalent result is that the Bayes factor $BF_{21} = m_2(\mathbf{X}_n)/m_1(\mathbf{X}_n)$ tends in probability $[P_1]$ to zero, where $\mathbf{X}_n = (X_1, \ldots, X_n)$.

The extension of this result to the case of a collection of models $\{P_i : i = 1, 2, \ldots\}$, for which the condition $\lim_{n \to \infty} m_i(\mathbf{X}_n)/m_1(\mathbf{X}_n) = 0$, $[P_1]$ holds for any $i \geq 2$, has been established by Dawid (1992). We note that this condition is satisfied when the model $P_1$ is nested in any other. For nonnested models the condition does not necessarily hold. As far as we know, a general consistency result for the Bayesian model selection procedure for nonnested models has not yet been established. This paper is a step forward in this direction and proves the consistency of Bayesian model selection procedures for normal linear models and a wide class of prior distributions, including the intrinsic priors.

For pairwise comparison between nested linear models the consistency of the intrinsic Bayesian procedure has already been established [Moreno and Girón (2005)]. The present paper is an extension of this result, and we prove here consistency of the intrinsic model posterior probabilities in the class of all linear models, where many of the models involved are nonnested. We also extend this result to a wide class of prior distributions. In proving consistency we take advantage of the nice asymptotic behavior of the Bayes factors arising from intrinsic priors. It is important to note we are assuming that the total number of regressors, $k$, is fixed and hence does not grow with $n$. For a consistency analysis where $k$ also grows with $n$, see Shao (1997).

The rest of the paper is organized as follows. In Section 2 we review methods for variable selection based on intrinsic priors and the expressions of Bayes factors and posterior model probabilities. In Section 3 we derive the sampling distributions of the statistic $\mathcal{B}_{ij}^n$, the statistic on which the Bayes factor for comparing two nested models depends, and we also describe its limiting behavior. This will be the tool we use in Section 4 to find out an asymptotic approximation of the Bayes factor for intrinsic priors, and to prove consistency of the variable selection procedure. Section 5 provides an evaluation of the intrinsic Bayes procedure and BIC for small sample sizes, and Section 6 contains a concluding discussion. There is also a short technical Appendix.

## 2. Intrinsic Bayesian procedures for variable selection.
Suppose that $Y$ represents an observable random variable and $X_1, X_2, \ldots, X_k$ a set of $k$ potential explanatory covariates related through the normal linear model

$$Y = \alpha_1 X_1 + \alpha_2 X_2 + \cdots + \alpha_k X_k + \varepsilon, \qquad \varepsilon \backsim \mathrm{N}(0, \sigma^2).$$

The variable selection problem consists of reducing the complexity of this model by identifying a subset of the $\alpha_i$ coefficients that have a zero value based on an available dataset $(\mathbf{y}, \mathbf{X})$, where $\mathbf{y}$ is a vector of observations of size $n$ and $\mathbf{X}$ an $n \times k$ design matrix of full rank.

This is by nature a model selection problem where we have to choose a model among the $2^k$ possible submodels of the above full one. It is common to set $X_1 = 1$ and $\alpha_1 \neq 0$ to include the intercept in any model. In this case the number of possible submodels is $2^{k-1}$. The class of models with $i$ regressors will be denoted as $\mathfrak{M}_i$ and hence the class of all possible submodels can be written as $\mathfrak{M} = \bigcup_i \mathfrak{M}_i$.

### 2.1. *Methods of encompassing.*
A fully Bayesian objective analysis for model comparison in linear regression has been given in Casella and Moreno

(2006). It consists of considering the pairwise model comparison between the full model $M_F$ and a generic submodel $M_i$[3] having $i$ ($< k$) nonzero regression coefficients. Formally, they test the hypothesis

$$(1) \qquad H_0 : \text{Model } M_i \quad \text{versus} \quad H_A : \text{Model } M_F.$$

Since $M_i$ is nested in the full model $M_F$, it is possible to derive the intrinsic priors for the parameters of both models. Then, in the space of models $\{M_i, M_F\}$ the intrinsic posterior probability of $M_i$ is computed using

$$P(M_i | \mathbf{y}, \mathbf{X}) = \frac{m_i(\mathbf{y}, \mathbf{X})}{m_i(\mathbf{y}, \mathbf{X}) + m_k(\mathbf{y}, \mathbf{X})} = \frac{BF_{ik}}{1 + BF_{ik}},$$

where $BF_{ik}$ is the Bayes factor for comparing model $M_i$ to model $M_F$. By doing this for every model an ordering of the set of models, in accordance to their posterior probabilities $\{P(M_i | \mathbf{y}, \mathbf{X}) = BF_{ik}/(1 + BF_{ik}), M_i \in \mathfrak{M}\}$, is obtained. The interpretation is that the submodel having the highest posterior probability is the most plausible reduction in complexity from the full model, the second highest the second-most plausible reduction and so on. This intrinsic Bayesian method for variable selection will be called *Variable Selection from Above* (VSA).

If we normalize the Bayes factors for intrinsic priors $\{BF_{ik}, i \geq 1\}$, we obtain a set of probabilities on the class $\mathfrak{M}$ as

$$(2) \qquad P(M_i; \mathbf{y}, \mathbf{X}) = \frac{BF_{ik}}{1 + \sum_{i' \neq k} BF_{i'k}}, \qquad M_i \in \mathfrak{M},$$

but we note that these probabilities are not true posterior probabilities of the models in the class $\mathfrak{M}$, although the ordering of the models they provide is exactly the same than that given by the above pairwise variable selection from above.

However, the manner of encompassing the models is not unique, and a quite natural alternative to VSA is to consider the pairwise model comparison between a generic submodel $M_j$ and the model

$$Y = \alpha_1 + \varepsilon, \qquad \varepsilon \curvearrowright \text{N}(\cdot | 0, \sigma^2),$$

that contains the intercept only, which is denoted as $M_1$. Formally, this comparison is made through the hypothesis test

$$(3) \qquad H_0 : \text{Model } M_1 \quad \text{versus} \quad H_A : \text{Model } M_j.$$

---

[3] We use $M_i$ to denote any model with $i$ regressors; there are $\binom{k-1}{i}$ such models. However, the development in the paper will be clear using this somewhat ambiguous, but simpler, notation.

Notice that $M_1$ is nested in $M_j$, for any $j$, so that the corresponding intrinsic priors can be derived. In the space of models $\{M_1, M_j\}$ the intrinsic posterior probability

$$P(M_j|\mathbf{y}, \mathbf{X}) = \frac{BF_{j1}}{1 + BF_{j1}}$$

is computed and it gives a new ordering of the models $\{M_j, M_j \in \mathfrak{M}\}$.

Although this alternative procedure is also based on multiple pairwise comparisons it is easy to see that it is equivalent to ordering the models according to the intrinsic model posterior probabilities computed in the space of all models $\mathfrak{M}$ as

(4)  $$P(M_j|\mathbf{y}, \mathbf{X}) = \frac{BF_{j1}}{1 + \sum_{j' \neq 1} BF_{j'1}}, \qquad M_j \in \mathfrak{M}.$$

This intrinsic Bayesian procedure will be called *Variable Selection from Below* (VSB), and has previously been considered by Girón, Moreno and Martinéz (2006).

For finite sample sizes, the orderings of the linear models provided by both VSA and VSB intrinsic Bayesian procedures are quite close to each other [Moreno and Girón (2008)].

2.2. *Intrinsic priors and Bayes factors.* The intrinsic priors utilized in the variable selection methods of Section 2.1 are defined from the comparison of two nested linear models, and we now give a general expression of the intrinsic priors and the Bayes factor associated with them.

Suppose we want to choose between the following two nested linear models

$$M_i : \mathbf{y} = \mathbf{X}_i \boldsymbol{\alpha}_i + \boldsymbol{\varepsilon}_i, \qquad \boldsymbol{\varepsilon}_i \sim \mathrm{N}_n(0, \sigma_i^2 \mathbf{I}_n)$$

and

$$M_j : \mathbf{y} = \mathbf{X}_j \boldsymbol{\beta}_j + \boldsymbol{\varepsilon}_j, \qquad \boldsymbol{\varepsilon}_j \sim \mathrm{N}_n(0, \sigma_j^2 \mathbf{I}_n).$$

We again can do this formally through the hypothesis test

(5)  $$H_0 : \text{Model } M_i \quad \text{versus} \quad H_A : \text{Model } M_j,$$

where $M_i$ is nested in $M_j$. Since the models are nested, this implies that the $n \times i$ design matrix $\mathbf{X}_i$ is a submatrix of the $n \times j$ design matrix $\mathbf{X}_j$, so that $\mathbf{X}_j = (\mathbf{X}_i | \mathbf{Z}_{ij})$. Then, model $M_j$ can be written as

$$M_j : \mathbf{y} = \mathbf{X}_i \boldsymbol{\beta}_i + \mathbf{Z}_{ij} \boldsymbol{\beta}_0 + \boldsymbol{\varepsilon}_j, \qquad \boldsymbol{\varepsilon}_j \sim \mathrm{N}_n(0, \sigma_j^2 \mathbf{I}_n).$$

Comparing model $M_i$ versus $M_j$ is equivalent to testing the hypothesis $H_0 : \beta_0 = 0$ against $H_1 : \beta_0 \neq 0$. A Bayesian setup for this problem is that of choosing between the Bayesian models

(6)
$$M_i : \mathrm{N}_n(\mathbf{y}|\mathbf{X}_i\boldsymbol{\alpha}_i, \sigma_i^2\mathbf{I}_n), \qquad \pi^N(\boldsymbol{\alpha}_i, \sigma_i) = \frac{c_i}{\sigma_i} \quad \text{and}$$

$$M_j : \mathrm{N}_n(\mathbf{y}|\mathbf{X}_j\boldsymbol{\beta}_j, \sigma_j^2\mathbf{I}_n), \qquad \pi^N(\boldsymbol{\beta}_j, \sigma_j) = \frac{c_j}{\sigma_j},$$

where $\pi^N$ denotes the improper reference prior and $c_i, c_j$ are arbitrary constants [Berger and Bernardo (1992)].

The direct use of improper priors for computing model posterior probabilities is not possible since they depend on the arbitrary constant $c_i/c_j$; however, they can be converted into suitable intrinsic priors [Berger and Pericchi (1996)]. Intrinsic priors for the parameters of the above nested linear models provide a Bayes factor [Moreno, Bertolino and Racugno (1998)] and, more importantly, posterior probabilities for the models $M_i$ and $M_j$, assuming that prior probabilities are assigned to them. Here we will use an objective assessment for this model prior probability, $P(M_i) = P(M_j) = 1/2$.

Application of the standard intrinsic prior methodology yields the intrinsic prior distribution for the parameters $\boldsymbol{\beta}_j, \sigma_j$ of model $M_j$, conditional on a fixed parameter point $\boldsymbol{\alpha}_i, \sigma_i$ of the reduced model $M_i$,

$$\pi^I(\boldsymbol{\beta}_j, \sigma_j|\boldsymbol{\alpha}_i, \sigma_i) = \frac{2}{\pi\sigma_i(1+\sigma_j^2/\sigma_i^2)}\mathrm{N}_j(\boldsymbol{\beta}_j|\tilde{\boldsymbol{\alpha}}_j, (\sigma_j^2+\sigma_i^2)\mathbf{W}_j^{-1}),$$

where $\tilde{\boldsymbol{\alpha}}_j' = (\mathbf{0}', \boldsymbol{\alpha}_i')$ with $\mathbf{0}$ being the null vector of $j-i$ components and

$$\mathbf{W}_j^{-1} = \frac{n}{j+1}(\mathbf{X}_j'\mathbf{X}_j)^{-1}.$$

The unconditional intrinsic prior for $(\boldsymbol{\beta}_j, \sigma_j)$ is obtained from $\pi^I(\boldsymbol{\beta}_j, \sigma_j) = \int \pi^I(\boldsymbol{\beta}_j, \sigma_j|\boldsymbol{\alpha}_i, \sigma_i)\pi^N(\boldsymbol{\alpha}_i, \sigma_i)\, d\boldsymbol{\alpha}_i\, d\sigma_i$, yielding the intrinsic priors for comparing models $M_i$ and $M_j$ as $\{\pi^N(\boldsymbol{\alpha}_i, \sigma_i), \pi^I(\boldsymbol{\beta}_j, \sigma_j)\}$. The computation of the Bayes factor to compare these models using the intrinsic priors is a straightforward calculation (see the Appendix) and turns out to be

(7)
$$BF_{ij}^n = \left(\frac{2}{\pi}(j+1)^{(j-i)/2}\right.$$
$$\left. \times \int_0^{\pi/2} \frac{\sin^{j-i}\varphi(n+(j+1)\sin^2\varphi)^{(n-j)/2}}{(n\mathcal{B}_{ij}^n + (j+1)\sin^2\varphi)^{(n-i)/2}}\, d\varphi\right)^{-1},$$

where the statistic $\mathcal{B}_{ij}^n$ is the ratio of the residual sum of squares

$$\mathcal{B}_{ij}^n = \frac{RSS_j}{RSS_i} = \frac{\mathbf{y}'(\mathbf{I}-\mathbf{H}_j)\mathbf{y}}{\mathbf{y}'(\mathbf{I}-\mathbf{H}_i)\mathbf{y}}.$$

Note that as $M_i$ is nested in $M_j$ the values of the statistic $\mathcal{B}_{ij}^n$ lie in the interval $[0, 1]$ and all of the above expressions are valid.

**3. Sampling distribution of $\mathcal{B}_{ij}^n$.** If we denote the true model by $M_T$, so that the distribution of the vector of observations $\mathbf{y}$ follows $\mathrm{N}_n(\mathbf{y}|\mathbf{X}_T\boldsymbol{\alpha}_T, \sigma_T^2\mathbf{I}_n)$, the sampling distribution of the statistic $\mathcal{B}_{ij}^n$ is given in the following theorem.

THEOREM 1.   *If $M_i$ is nested in $M_j$ and $M_T$ is the true model, then the sampling distribution of $\mathcal{B}_{ij}^n$ is the doubly noncentral beta distribution*

$$\mathcal{B}_{ij}^n|M_T \backsim \mathrm{Be}\left(\frac{n-j}{2}, \frac{j-i}{2}; \lambda_1, \lambda_2\right),$$

*where the noncentrality parameters are*

$$\lambda_1 = \frac{1}{2\sigma_T^2}\boldsymbol{\alpha}_T'\mathbf{X}_T'(\mathbf{I} - \mathbf{H}_j)\mathbf{X}_T\boldsymbol{\alpha}_T$$

*and*

$$\lambda_2 = \frac{1}{2\sigma_T^2}\boldsymbol{\alpha}_T'\mathbf{X}_T'(\mathbf{H}_j - \mathbf{H}_i)\mathbf{X}_T\boldsymbol{\alpha}_T.$$

PROOF.   The quadratic form of the denominator of the $\mathcal{B}_{ij}^n$ can be decomposed as

$$\mathbf{y}'(\mathbf{I} - \mathbf{H}_i)\mathbf{y} = \mathbf{y}'(\mathbf{I} - \mathbf{H}_j)\mathbf{y} + \mathbf{y}'(\mathbf{H}_j - \mathbf{H}_i)\mathbf{y}.$$

As the matrices $(\mathbf{I} - \mathbf{H}_i)$ and $(\mathbf{H}_j - \mathbf{H}_i)$ are idempotent of ranks $n-j$ and $j-i$, respectively, it follows from the generalized Cochran theorem that the quadratic form $\mathbf{y}'(\mathbf{I} - \mathbf{H}_j)\mathbf{y}$ and $\mathbf{y}'(\mathbf{H}_j - \mathbf{H}_i)\mathbf{y}$ are independent and distributed as $\chi'^2(n-j; \lambda_1)$ and $\chi'^2(j-i; \lambda_2)$, respectively. From this the distribution of the statistic $\mathcal{B}_{ij}^n$ follows, and Theorem 1 is proved.   □

Note that the models $M_i$ and $M_j$ need not be nested in the true model $M_T$, and the true model is not necessarily nested in $M_i$ or $M_j$. However, the distribution of $\mathcal{B}_{ij}^n$ simplifies whenever $M_i$ or $M_j$ is the true model. Thus we have the following corollary.

COROLLARY 1.   (i) *If the smaller model $M_i$ is the true one, then*

$$\mathcal{B}_{ij}^n|M_i \backsim \mathrm{Be}\left(\frac{n-j}{2}, \frac{j-i}{2}\right).$$

(ii) *If the larger model $M_j$ is the true one, then*

$$\mathcal{B}_{ij}^n|M_j \backsim \mathrm{Be}\left(\frac{n-j}{2}, \frac{j-i}{2}; 0, \lambda\right),$$

*where*

$$\lambda = \frac{1}{2\sigma_j^2}\boldsymbol{\alpha}_j'\mathbf{X}_j'(\mathbf{H}_j - \mathbf{H}_i)\mathbf{X}_j\boldsymbol{\alpha}_j.$$

PROOF. Part (i) follows from the fact that $\mathbf{X}_i'\mathbf{H}_j = \mathbf{X}_i'\mathbf{H}_i$ and part (ii) from $\mathbf{X}_j'(\mathbf{H}_j - \mathbf{H}_i) = \mathbf{X}_j'(\mathbf{I} - \mathbf{H}_i)$.  □

The limiting value of $\mathcal{B}_{ij}^n$ is important because it bears directly on the evaluation of the consistency of the Bayes factors. That value is given in the following theorem.

THEOREM 2. *Let $\{X_n, n \geq 1\}$ be a sequence of random variables with distribution* $\mathrm{Be}((n - \alpha_0)/2, \beta_0/2; n\delta_1, n\delta_2)$, *where $\alpha_0, \beta_0, \delta_1, \delta_2$ are positive constants. Then:*

(i) *the sequence $X_n$ converges in probability to the constant*

$$\frac{1 + \delta_1}{1 + \delta_1 + \delta_2};$$

(ii) *if $\delta_1 = \delta_2 = 0$, then $X_n$ degenerates in probability to 1. However, the random variable $-n/2 \log X_n$ does not degenerate and has an asymptotic Gamma distribution,* $\mathrm{Ga}(\beta_0, 1)$.

PROOF. Part (i). By definition $X_n$ is

$$X_n = \left(1 + \frac{\chi'^2(\beta_0, n\delta_2)}{\chi'^2(n - \alpha_0, n\delta_1)}\right)^{-1},$$

where $\chi'^2(\beta_0, n\delta_2)$ and $\chi'^2(n - \alpha_0, n\delta_1)$ are independent random variables with noncentral chi-square distributions. If we divide the numerator and denominator by $n$ we get

$$X_n = \left(1 + \frac{V_n}{W_n}\right)^{-1},$$

where $V_n = \chi'^2(\beta_0, n\delta_2)/n$ and $W_n = \chi'^2(n - \alpha_0, n\delta_1)/n$. Their means and variances are

$$E(V_n) = \delta_2 + \frac{\beta_0}{n}, \qquad E(W_n) = 1 + \delta_1 - \frac{\alpha_0}{n}$$

and

$$\mathrm{Var}(V_n) = \frac{4\delta_2}{n} + \frac{2\beta_0}{n^2}, \qquad \mathrm{Var}(W_n) = \frac{2(1 + \delta_1)}{n} - \frac{2\alpha_0}{n^2}.$$

Since the variances go to zero as $n$ goes to infinity, $X_n$ degenerates in probability to $(1 + \delta_1)/(1 + \delta_1 + \delta_2)$ as asserted.

The remainder of the proof is straightforward and hence is omitted.  □

**4. Consistency of the VSA and VSB intrinsic Bayesian procedures.** The steps in proving consistency of the intrinsic Bayesian procedures are:

1. Derive an asymptotic approximation for the Bayes factor for nested models given in expression (7).

2. From this approximation derive another that is valid for any arbitrary pair of models.

3. Use Theorems 1 and 2 to prove consistency of the VSB procedure.

It will also be seen that the asymptotic behavior of the Bayes factor for VSA is exactly the same as VSB, and hence the consistency of the former procedure also holds.

This is a useful property of the intrinsic methodology for variable selection since any way of encompassing the models to derive the intrinsic priors produces essentially the same answer for finite sample sizes and for large sample sizes.

4.1. *Asymptotic approximation of $BF_{ij}^n$.* For large $n$, we can get an approximation of $BF_{ij}^n$ of (7) that is valid whenever model $M_i$ is nested in $M_j$. The approximation turns out to be equivalent to the Schwarz (1978) Bayes factor approximation.

THEOREM 3. *When $M_i$ is nested in $M_j$, for large values of $n$ the Bayes factor given in (7) can be approximated by*

$$(8) \quad BF_{ij}^n \approx \frac{\pi}{2}(j+1)^{(i-j)/2} I(\mathcal{B}_{ij}^n)^{-1} \exp\left(\frac{j-i}{2}\log n + \frac{n-i}{2}\log \mathcal{B}_{ij}^n\right),$$

*where*

$$I(\mathcal{B}_{ij}^n) = \int_0^{\pi/2} \sin^{j-i}(\varphi) \exp\left[\frac{j+1}{2}\sin^2(\varphi)\left(1 - \frac{1}{\mathcal{B}_{ij}^n}\right)\right] d\varphi$$

$$= \frac{1}{2}\operatorname{Be}\left(\frac{1}{2}, \frac{j-i+1}{2}\right)$$

$$\times {}_1F_1\left(\frac{j-i+1}{2}; \frac{j-i+2}{2}; \frac{j+1}{2}\left(1 - \frac{1}{\mathcal{B}_{ij}^n}\right)\right),$$

*and ${}_1F_1(a;b;z)$ denotes the Kummer confluent hypergeometric function [see Abramowitz and Stegun (1972), Chapter 13].*

PROOF. We can write the integrand of (7) as

$$\sin^{j-i}\varphi \exp\left\{\frac{n-j}{2}\left[\log n + \log\left(1 + \frac{j+1}{n}\sin^2\varphi\right)\right]\right\}$$

$$\times \exp\left\{\frac{i-n}{2}\left[\log n + \log \mathcal{B}_{ij}^n + \log\left(1 + \frac{j+1}{n\mathcal{B}_{ij}^n}\sin^2\varphi\right)\right]\right\}$$

$$= \sin^{j-i}\varphi \exp\left(\frac{i-j}{2}\log n + \frac{i-n}{2}\log \mathcal{B}_{ij}^n\right)$$

$$\times \frac{(1+(j+1)/n\sin^2\varphi)^{(n-j)/2}}{(1+(j+1)/(n\mathcal{B}_{ij}^n)\sin^2\varphi)^{(n-i)/2}}.$$

For large $n$ the numerator of the last factor can be approximated by

$$\left(1+\frac{j+1}{n}\sin^2\varphi\right)^{(n-j)/2} \approx \exp\left\{\frac{j+1}{2}\sin^2\varphi\right\},$$

and the denominator by

$$\left(1+\frac{j+1}{n\mathcal{B}_{ij}^n}\sin^2\varphi\right)^{(n-i)/2} \approx \exp\left\{\frac{j+1}{2\mathcal{B}_{ij}^n}\sin^2\varphi\right\}.$$

Therefore, for large $n$ the integrand can be approximated by

$$\sin^{j-i}\varphi \exp\left(\frac{i-j}{2}\log n + \frac{i-n}{2}\log \mathcal{B}_{ij}^n\right)\exp\left(\frac{j+1}{2}\sin^2\varphi\left(1-\frac{1}{\mathcal{B}_{ij}^n}\right)\right),$$

and thus the Bayes factor (7) by

$$BF_{ij}^n \approx \frac{\pi}{2}(j+1)^{(i-j)/2}I(\mathcal{B}_{ij}^n)^{-1}\exp\left(\frac{j-i}{2}\log n + \frac{n-i}{2}\log \mathcal{B}_{ij}^n\right),$$

where

$$I(\mathcal{B}_{ij}^n) = \int_0^{\pi/2}\sin^{j-i}\varphi \exp\left[\frac{j+1}{2}\sin^2\varphi\left(1-\frac{1}{\mathcal{B}_{ij}^n}\right)\right]d\varphi.$$

This proves Theorem 3. □

We note that $I(\mathcal{B}_{ij}^n)^{-1}$ has a finite value for all values of the statistic $\mathcal{B}_{ij}^n$ except when it goes to zero. However, we can see in the proof of Theorem 4 that $\mathcal{B}_{ij}^n$ tends to a strictly positive number with probability 1 as $n \to \infty$ [see expression (14)], so $I(\mathcal{B}_{ij}^n)^{-1}$ is finite for all $n$.

Therefore, $BF_{ij}^n$ can be approximated, up to a multiplicative constant, by the exponential function in (8). This exponential function turns out to be the Schwarz approximation $S_{ij}^n$ to the Bayes factor for comparing linear models [Schwarz (1978)]. Of course, the normal linear models are regular so the Laplace approximation can be applied to obtain the Schwarz approximation although for intrinsic priors the ratio $BF_{ij}^n/S_{ij}^n$ does not go to 1 [this holds only for particular priors; see Kass and Wasserman (1995)].

However, for proving consistency we can ignore terms of constant order and the Bayes factor for intrinsic priors can be approximated by the Schwarz approximation

$$(9) \qquad BF_{ij}^n \approx S_{ij}^n = \exp\left(\frac{j-i}{2}\log n + \frac{n}{2}\log \mathcal{B}_{ij}^n\right).$$

We note that $S_{ij}^n$ could provide a crude approximation to $BF_{ij}^n$ for small or moderate sample sizes. In Section 5 we look at small-sample behavior of both the Schwarz approximation and the Bayes factor for intrinsic priors.

4.2. *Consistency of the VSB intrinsic Bayesian procedure.* Given an arbitrary model $M_j$ and the true model $M_T$ in the class $\mathfrak{M}_T$, we will assume the design matrix of the linear models satisfy the following condition (D): the matrix

(10) $$\mathbf{S}_{jT} = \lim_{n \to \infty} \frac{\mathbf{X}_T'(\mathbf{I} - \mathbf{H}_j)\mathbf{X}_T}{n}$$

is a positive semidefinite matrix. This is not a too demanding condition as the following example shows.

EXAMPLE 1 [Berger and Pericchi (2004)]. Consider the case of testing whether the slope of a linear regression is zero. Suppose that the true model $M_T$ is the model with regression coefficients $(\alpha_1, \alpha_2)$, and thus there is only one alternative model $M_1$, the model with only the intercept term $\alpha_1$. Suppose that there are $2n + 1$ observations yielding the design matrix

$$\mathbf{X}^t = \begin{pmatrix} 1 & \dots & 1 & 1 & \dots & 1 & 1 \\ 0 & \dots & 0 & \delta & \dots & \delta & 1 \end{pmatrix},$$

where $\delta$ is different from zero. Easy calculations show that

$$\mathbf{S}_{1T} = \lim_{n \to \infty} \frac{\mathbf{X}_T'(\mathbf{I} - \mathbf{H}_1)\mathbf{X}_T}{2n+1} = \begin{pmatrix} 0 & 0 \\ 0 & \delta^2/4 \end{pmatrix},$$

which obviously is a positive semidefinite matrix for any positive $|\delta|$, no matter how close to zero it is.

Thus, condition (D) is satisfied even when the samples are coming from a model $M_T$, which is as close to $M_1$ as we want.

To characterize the asymptotic behavior of the model posterior probabilities, we can work with $BF_{ij}^n$ of (8), ignoring the positive terms that do not depend on $n$ as we are only interested in limiting values of 0 or $\infty$.

To test the hypothesis (3) with data $(\mathbf{y}, \mathbf{X})$, we note that the intrinsic model posterior probability of model $M_j$, defined in the class of all models $\mathfrak{M}$ given by (4), is an increasing function of $BF_{j1}$, where $BF_{j1}$ denotes the Bayes factor for intrinsic priors for comparing the nested models $M_1$ versus $M_j$. Hence, from the asymptotic approximation (8) we can write

$$P(M_j|\mathbf{y}, \mathbf{X}) = \frac{BF_{j1}}{1 + \sum_{j' \neq 1} BF_{j'1}}$$

$$= \left( c_{j1} I(\mathcal{B}_{1j}^n)^{-1} \exp\left\{ -\frac{j-1}{2} \log n - (n/2) \log \mathcal{B}_{1j}^n \right\} \right)$$

(11)
$$\times \left(1 + \sum_{j' \neq 1} c_{j'1} I(\mathcal{B}_{1j'}^n)^{-1} \exp\left\{-\frac{j'-1}{2}\log n \right.\right.$$
$$\left.\left. - (n/2)\log \mathcal{B}_{1j'}^n\right\}\right)^{-1}.$$

Similarly, for the true model $M_T$ we can write

$$P(M_T|\mathbf{y}, \mathbf{X})$$
$$= \frac{c_{T1} I(\mathcal{B}_{1T}^n)^{-1} \exp\{-((T-1)/2)\log n - (n/2)\log \mathcal{B}_{1T}^n\}}{1 + \sum_{j' \neq 1} c_{j'1} I(\mathcal{B}_{1j'}^n)^{-1} \exp\{-((j'-1)/2)\log n - (n/2)\log \mathcal{B}_{1j'}^n\}},$$

where $c_{j1}$ and $c_{T1}$ do not depend on $n$, and $I(\mathcal{B}_{1j}^n)^{-1}$ and $I(\mathcal{B}_{1T}^n)^{-1}$ are finite for all $n$. We are concerned with the limiting behavior of the ratio of these two probabilities, and specifically if the limit is 0 or $\infty$. Thus, in the following we can ignore the finite terms and approximate the ratio with

(12)
$$\frac{P(M_j|\mathbf{y}, \mathbf{X})}{P(M_T|\mathbf{y}, \mathbf{X})} \approx \exp\left\{\frac{T-j}{2}\log n + \frac{n}{2}\log\frac{\mathcal{B}_{1T}^n}{\mathcal{B}_{1j}^n}\right\},$$

because the denominators cancel. (As a curiosity, note that this formula provides an exact approximation to the ratio for the case when $M_j = M_T$, when its value is exactly equal to one.)

We now have the following theorem.

THEOREM 4. *In the class of linear models $\mathfrak{M}$ with design matrices satisfying condition (D), the intrinsic Bayesian variable selection procedure VSB is consistent. That is, when sampling from $M_T$ we have that*

$$\frac{P(M_j|\mathbf{y}, \mathbf{X})}{P(M_T|\mathbf{y}, \mathbf{X})} \to 0, \qquad [P_t],$$

*whenever the model $M_j \neq M_T$.*

PROOF. Assuming $M_T \neq M_1$, from Corollary 1(ii), we have that

$$\mathcal{B}_{1T}^n|M_T \backsim \mathrm{Be}\left(\frac{n-T}{2}, \frac{T-1}{2}; 0, \lambda\right),$$

where

$$\lambda = \frac{1}{2\sigma_T^2}\boldsymbol{\alpha}_T' \mathbf{X}_T'(\mathbf{I} - \mathbf{H}_1)\mathbf{X}_T\boldsymbol{\alpha}_T$$

and from Theorem 1 that

$$\mathcal{B}_{1j}^n|M_T \backsim \mathrm{Be}\left(\frac{n-j}{2}, \frac{j-1}{2}; \lambda_1, \lambda_2\right),$$

where the noncentrality parameters are

$$
\begin{aligned}
\lambda_1 &= \frac{1}{2\sigma_T^2}\boldsymbol{\alpha}_T'\mathbf{X}_T'(\mathbf{I}-\mathbf{H}_j)\mathbf{X}_T\boldsymbol{\alpha}_T, \\
\lambda_2 &= \frac{1}{2\sigma_T^2}\boldsymbol{\alpha}_T'\mathbf{X}_T'(\mathbf{H}_j-\mathbf{H}_1)\mathbf{X}_T\boldsymbol{\alpha}_T.
\end{aligned}
\tag{13}
$$

From Theorem 2(i), we have

$$
\begin{aligned}
\mathcal{B}_{1T}^n|M_T &\to \frac{1}{1+1/(2\sigma_T^2)\boldsymbol{\alpha}_T'\mathbf{S}_{1T}\boldsymbol{\alpha}_T} \quad \text{and} \\
\mathcal{B}_{1j}^n|M_T &\to \frac{1+1/(2\sigma_T^2)\boldsymbol{\alpha}_T'\mathbf{S}_{jT}\boldsymbol{\alpha}_T}{1+1/(2\sigma_T^2)\boldsymbol{\alpha}_T'\mathbf{S}_{1T}\boldsymbol{\alpha}_T},
\end{aligned}
\tag{14}
$$

so that

$$
\frac{\mathcal{B}_{1T}^n}{\mathcal{B}_{1j}^n}\Big|M_T \to \frac{1}{1+1/(2\sigma_T^2)\boldsymbol{\alpha}_T'\mathbf{S}_{jT}\boldsymbol{\alpha}_T} < 1.
$$

Therefore, the expression

$$
\frac{n}{2}\log\frac{\mathcal{B}_{1T}^n}{\mathcal{B}_{1j}^n}
$$

goes to $-\infty$ with order $O(n)$. This means that expression (12) converges to zero regardless of whether $T-j$ is positive or negative.

When $M_T = M_1$, then for any $j > 1$ we have

$$
P(M_j|\mathbf{y},\mathbf{X}) \propto BF_{j1}^n \approx \exp\left(-\frac{j-1}{2}\log n - \frac{n}{2}\log\mathcal{B}_{1j}^n\right).
$$

From Corollary 1(i) and Theorem 2(ii) it follows that $-n/2\log\mathcal{B}_{1j}^n$ is asymptotically distributed as a Gamma distribution. Therefore, for any $j > 1$, $P(M_j|\mathbf{y},\mathbf{X})$ tends, in probability, to zero. The proof is complete. $\square$

4.3. *Consistency of the VSA intrinsic Bayesian procedure.* In the VSA intrinsic Bayesian procedure we use the fact that every model $M_j$ is nested in the full model $M_k$. Then, for large values of $n$ the posterior probability of model $M_j$ in the space of models $\{M_j, M_k\}$ is proportional to

$$
P(M_j|\mathbf{y},\mathbf{X}) \propto BF_{jk}^n \approx \exp\left(\frac{k-j}{2}\log n + \frac{n}{2}\log\mathcal{B}_{jk}^n\right).
$$

Similarly, for the true model $M_T$ we have

$$
P(M_T|\mathbf{y},\mathbf{X}) \propto BF_{Tk}^n \approx \exp\left(\frac{k-T}{2}\log n + \frac{n}{2}\log\mathcal{B}_{Tk}^n\right).
$$

Thus, the ratio of Bayes factors can be approximated by

$$\frac{P(M_j|\mathbf{y}, \mathbf{X})}{P(M_T|\mathbf{y}, \mathbf{X})} \propto \frac{BF_{jk}^n}{BF_{Tk}^n} \approx \exp\left(\frac{T-j}{2}\log n + \frac{n}{2}\log\frac{\mathcal{B}_{1T}^n}{\mathcal{B}_{1j}^n}\right)$$

where the last expression is exactly that given in (12) so that it tends to zero for any $j \geq 1$. We thus have the following corollary to Theorem 4.

COROLLARY 2. *In the class of linear models $\mathfrak{M}$ with design matrices satisfying condition (D), the intrinsic Bayesian variable selection procedure VSA is consistent. That is, when sampling from $M_T$ we have that*

$$\frac{P(M_j|\mathbf{y}, \mathbf{X})}{P(M_T|\mathbf{y}, \mathbf{X})} \to 0, \qquad [P_t],$$

*whenever the model $M_j \neq M_T$.*

Recall that in Section 2.1 we noted that for VSA, the probabilities

$$P(M_i|\mathbf{y}, \mathbf{X}) = \frac{BF_{ik}^n}{1 + \sum_{i' \neq k} BF_{i'k}^n}, \qquad M_i \in \mathfrak{M},$$

were not true posterior probabilities of the models in the class $\mathfrak{M}$. However, from Corollary 2, this set of probabilities [utilized as a tool for variable selection in Casella and Moreno (2006)], is a consistent sequence of probabilities. Further, we recall that the ordering of the models they provide is exactly the same as that given by the VSA pairwise variable selection. Therefore, the intrinsic model posterior probabilities from above form a set of consistent probabilities in the class of all linear models $\mathfrak{M}$.

4.4. *Extensions.* The consistency of the intrinsic Bayesian variable selection procedure for the class of linear models can be extended to any other Bayesian procedure for a wide class of prior distributions. We observe that all we have used to prove consistency of the intrinsic Bayesian procedures is the Schwarz approximation, and the distribution of the ratio of the residuals of two nested linear models when sampling from a linear model that does not necessarily coincide with any of the two. Therefore, for any prior for which the Schwarz approximation for linear models is valid, the consistency of the associated Bayesian procedure can be asserted. Hence, we can prove the following theorem.

THEOREM 5. *In the class of linear models $\mathfrak{M}$ with design matrices satisfying condition (D), assume that the priors $\pi_i$, $\pi_j$ for any $i, j$, are such that*

$$0 < \lim_{n \to \infty} \frac{\pi_i(\hat{\alpha}_i, \hat{\sigma}_i)}{\pi_j(\hat{\alpha}_j, \hat{\sigma}_j)} < \infty, \qquad [P_T],$$

*where $\hat{\alpha}_i, \hat{\sigma}_i$ and $\hat{\alpha}_j, \hat{\sigma}_j$ are the respective MLEs. Then the Bayesian variable selection procedure is consistent, that is, when sampling from $M_T \in \mathfrak{M}$, we have that*

$$\frac{P(M_j|\mathbf{y}, \mathbf{X})}{P(M_T|\mathbf{y}, \mathbf{X})} \to 0, \qquad [P_t],$$

*whenever the model $M_j \neq M_T$.*

We note that priors of the form $\pi_i^N(\boldsymbol{\alpha}_i, \sigma_i^q) = c_i/\sigma_i^q$, where $q$ is a positive number, which includes the reference priors for $q = 1$ and the Jeffreys priors for $q = i$, satisfy the condition required in Theorem 5. Indeed, from (14), it follows that

$$\lim_{n\to\infty} \frac{\pi_i^N(\hat{\boldsymbol{\alpha}}_i, \hat{\sigma}_i)}{\pi_j^N(\hat{\boldsymbol{\alpha}}_j, \hat{\sigma}_j)} = \left(\frac{c_i}{c_j} \lim_{n\to\infty} \mathcal{B}_{ij}^n\right)^{q/2}$$

$$= \left(\frac{c_i}{c_j}\right)^{q/2} \left(\frac{2\sigma_T^2 + \boldsymbol{\alpha}_T' \mathbf{S}_{jT} \boldsymbol{\alpha}_T}{2\sigma_T^2 + \boldsymbol{\alpha}_T' \mathbf{S}_{iT} \boldsymbol{\alpha}_T}\right)^{q/2}, \qquad [P_T]$$

which clearly is a real positive quantity.

Hence, even though for finite sample sizes the above priors only provide Bayes factors defined up to a multiplicative constant, asymptotically they behave consistently.

**5. Small sample comparisons.** Although for large sample sizes the variable selection procedure based on the Bayes factor for intrinsic priors is equivalent to that based on the Schwarz approximation, an open question is how good the Schwarz asymptotic approximation and the Bayes factor for intrinsic priors behave for small or moderate sample sizes. To answer this question we recall that, in the case of encompassing from below, the ordering of the models provided by the pairwise intrinsic model posterior probabilities

$$P(M_j|\mathbf{y}, \mathbf{X}) = \frac{B_{j1}^n}{1 + B_{j1}^n} \qquad \text{for } j \geq 2$$

is exactly the same as that provided by the intrinsic model posterior probabilities in the whole space $\mathfrak{M}$.

Therefore, for comparing the intrinsic Bayes factor $B_{ij}^n$ and the Schwarz approximation $S_{ij}^n$ for any $i$ and $j$ it is enough to compare $B_{j1}^n$ and $S_{j1}^n$ for $j \geq 2$. It seems appropriate to compare $B_{j1}^n$ and $S_{j1}^n$ in a probabilistic scale, that is, to compare the intrinsic posterior model probability $P(M_j|\mathbf{y}, \mathbf{X})$ and the Schwarz approximation posterior probability

$$P^S(M_j|\mathbf{y}, \mathbf{X}) = \frac{S_{j1}^n}{1 + S_{j1}^n} \qquad \text{for } j \geq 2.$$

TABLE 1
*Type* I *error probabilities for the intrinsic procedure and the Schwarz approximation. In each cell, the left probability is the Type* I *error of the intrinsic procedure and the right probability is the Type* I *error of the Schwarz approximation*

|           | $n = 7$    | $n = 10$   | $n = 15$      | $n = 40$     | $n = 80$       |
|-----------|------------|------------|---------------|--------------|----------------|
| $j = 2$   | 0.16, 0.26 | 0.13, 0.19 | 0.10, 0.130   | 0.06, 0.06   | 0.04, 0.04     |
| $j = 3$   | 0.19, 0.33 | 0.14, 0.20 | 0.099, 0.114  | 0.04, 0.03   | 0.02, 0.02     |
| $j = 4$   | 0.23, 0.42 | 0.16, 0.22 | 0.104, 0.102  | 0.03, 0.02   | 0.02, 0.02     |
| $j = 5$   | 0.29, 0.55 | 0.18, 0.25 | 0.111, 0.097  | 0.03, 0.01   | 0.01, 0.002    |
| $j = 6$   | 0.40, 0.75 | 0.21, 0.31 | 0.121, 0.097  | 0.03, 0.006  | 0.01, 0.001    |
| $j = 9$   |            | 0.41, 0.71 | 0.17, 0.15    | 0.03, 0.002  | 0.007, $\simeq 0$ |
| $j = 12$  |            |            | 0.26, 0.36    | 0.04, 0.001  | $\simeq 0$, $\simeq 0$ |
| $j = 38$  |            |            |               | 0.32, 0.46   | 0.001, $\simeq 0$ |
| $j = 78$  |            |            |               |              | 0.32, 0.44     |

A model selection procedure operates by choosing the model with the highest value of the criterion, so in our case this is equivalent to accepting model $M_j$, and hence rejecting $M_1$, when the posterior probability of $M_j$ is greater than $1/2$. It is important to realize that this is not the way a classical frequentist hypothesis test is set up. In the classical case a test is calibrated to a specified Type I error $\alpha$, and then the power is examined. The model selector is defined by the decision rule, and for the given rule we can examine the resulting Type I and II errors to assess how the model selector is controlling them.

We recall that both the intrinsic posterior probability $P(M_j|\mathbf{y}, \mathbf{X})$ and the Schwarz approximation $P^S(M_j|\mathbf{y}, \mathbf{X})$ depend on the sample observations $(\mathbf{y}, \mathbf{X})$ through the statistic $\mathcal{B}_{1j}^n$. Therefore, any point in the regions

$$R_{1j}(n) = \{\mathcal{B}_{1j}^n : P(M_j|\mathcal{B}_{1j}^n) \geq 1/2\} \qquad \text{and}$$

$$R_{1j}^S(n) = \{\mathcal{B}_{1j}^n : P^S(M_j|\mathcal{B}_{1j}^n) \geq 1/2\}$$

contain empirical evidence in favor of $M_j$ under the intrinsic Bayesian procedure and the Schwarz approximation.

Since $M_1$ is nested in $M_j$ for any $j \geq 2$, it follows that $R_{1j}(n) \subset (0, 1)$, and $R_{1j}^S(n) \subset (0, 1)$. Furthermore, $R_{1j}(n)$ and $R_{2j}(n)$ are intervals since $P(M_j|\mathcal{B}_{1j}^n)$ and $P^S(M_j|\mathcal{B}_{1j}^n)$ are monotone increasing functions of $\mathcal{B}_{1j}^n$.

The distribution of $\mathcal{B}_{1j}^n$ is easily computed (see Corollary 1), and we can examine the Type I errors of the intrinsic Bayesian variable selection procedure and the Schwarz approximation, respectively. For a range of values of $j$ and sample sizes $n > j$, Table 1 presents the Type I error probabilities under the intrinsic Bayesian procedure and the Schwarz approximation.

We see in Table 1 that for small sample sizes the Schwarz approximation has a very high Type I error (as high as 75%), which soon becomes very small as $n$ increases. Thus, the Schwarz approximation will be biased away from the null model for small $n$, or more generally, in the cases where $j$ is close to $n$. As $n$ increases the Type I error goes rapidly to 0, and the Schwarz approximation will then be biased toward the null model. In contrast, the intrinsic procedure has a less variable Type I error, being smaller than that of the Schwarz approximation for small $n$ and somewhat larger for large $n$.

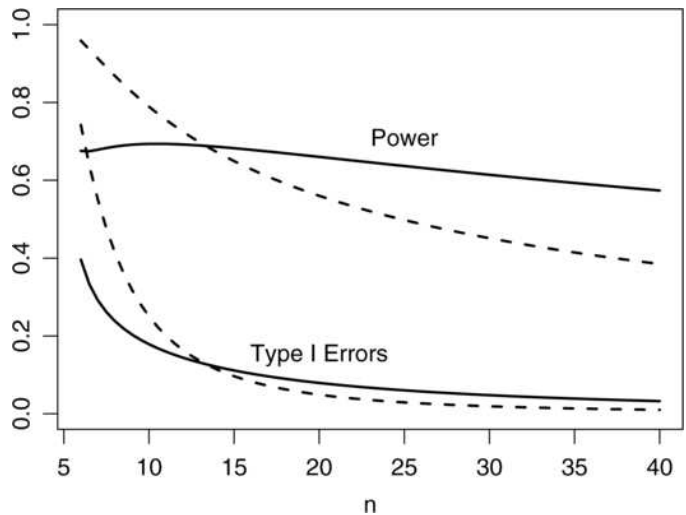

Fig. 1. *For $j = 5$ and $n = 6, \ldots, 40$, Type I errors and power curves of the intrinsic procedure (solid) and Schwarz approximation (dashed) as a function of n. The power curves are computed for noncentrality parameter $\lambda = 10$.*

Examination of Figure 1 shows a very interesting story. There, we plotted Type I errors and power as a function of $n$ for $j = 5$, which was chosen as a representative case. Note that the decrease in the power, as a function of $n$, reflects the fact that the Type I error decreases as a function of $n$.

For small $n$ the Schwarz approximation has higher power resulting from its large Type I error, while the intrinsic procedure tends to moderate both errors. As $n$ increases, both Type I errors decrease, with the more dramatic decrease being that of the Schwarz approximation. The Type I errors cross at $n = 13$, and for $n > 13$ the intrinsic procedure has higher power, reaching 0.573 at $n = 40$ versus 0.385 for the Schwarz approximation. The interesting point is that, although the intrinsic procedure has higher Type I error, both Type I errors are very small (e.g., at $n = 29$ they are 0.05 and 0.02). However, the effect of Schwarz approximation, by driving the Type I error so close to zero, is a dramatic decrease in power. Thus, the intrinsic procedure does a much better job of controlling the errors. By moderating the Type I error it avoids the faults of the Schwarz approximation, which has very large Type I error for small $n$, and for large $n$ decreases the Type I error to an unnecessarily low value to the detriment of its power.

**6. Discussion.** It has long been known that when choosing between two models, when one of which is true, selecting according to Bayes factors provides a consistent decision function in the sense that the *frequentist* probability of selecting the true model approaches 1 as $n \to \infty$. In this paper, for the case of variable selection, we have extended this result to selection among an entire class of linear models and a wide class of priors, and shown that selecting according to Bayes factors yields a decision rule with the property that the frequentist probability of selecting the true model approaches 1 as $n \to \infty$, and the frequentist probability of selecting any other model approaches 0 as $n \to \infty$.

We have, specifically, worked with intrinsic priors, although our results hold for a wide class of priors. However, intrinsic priors provide a type of objective Bayesian prior for the testing problem. They seem to be among the most diffuse priors that are possible to use in testing, without encountering problems with indeterminate Bayes factors, which was the original impetus for the development of Berger and Pericchi ([1996](#)). Moreover, they do not suffer from "Lindley paradox" behavior. Thus, we believe they are a very reasonable choice for experimenters looking for an objective Bayesian analysis with a frequentist guarantee. This is very much in the spirit of the *calibrated Bayesian*, as described by Little ([2006](#)).

Intrinsic priors have been used successfully in both variable selection and changepoint problems [Casella and Moreno ([2006](#)), Girón, Moreno and Martínez ([2006](#)), Girón, Moreno and Casella ([2007](#))], where excellent small-sample properties were exhibited. Some other properties of the variable selection rules considered here are as follows:

1. All models $M_j$ that contain model $M_T$, and hence have $\lambda_1 = 0$ [see ([13](#))], will have the same value of $\mathcal{B}_{1T}^n | M_T$ in ([14](#)). This means that the posterior probability of models $M_j$ that contain model $M_T$ ([11](#)) is decreasing in $j$, and models with larger $j$ will have smaller probabilities. Thus, VSB will tend to select smaller models. The same holds for VSA.

2. To gain further insight in the large-sample approximation of the Bayes factors for comparing arbitrary models, say $M_j$ and $M_{j'}$, we look a bit closer at the importance of some geometric considerations in the space of all models, as the one played by a distance that we can define between a generic model $M_j$ and the true, though unknown, model $M_T$.

If we define this distance as

$$\delta(M_j, M_T) = \frac{\alpha_T' \mathbf{S}_{jT} \alpha_T}{\sigma_T^2},$$

we note that it is equal to 0 if either $M_j = M_T$ or $M_T$ is nested in $M_j$; otherwise, it is strictly positive by condition (D). Also, if model $M_i$ is nested in $M_j$ then $\delta(M_i, M_T) < \delta(M_j, M_T)$, because $\mathbf{H}_j - \mathbf{H}_i$ is positive semidefinite.

3. From (11) we have that

$$\frac{P(M_j|\mathbf{y}, \mathbf{X})}{P(M_{j'}|\mathbf{y}, \mathbf{X})} \approx \exp\left(\frac{j'-j}{2}\log n - \frac{n}{2}\log\frac{\mathcal{B}_{1j}^n}{\mathcal{B}_{1j'}^n}\right)$$

and from (14)

$$\log\frac{\mathcal{B}_{1j}^n}{\mathcal{B}_{1j'}^n}\Big|M_T \to \log\frac{1+\delta(M_j, M_T)/2}{1+\delta(M_{j'}, M_T)/2}.$$

Hence,

$$\frac{P(M_j|\mathbf{y}, \mathbf{X})}{P(M_{j'}|\mathbf{y}, \mathbf{X})}\Big|M_T \approx \exp\left(\frac{j'-j}{2}\log n - \frac{n}{2}\log\frac{1+\delta(M_j, M_T)/2}{1+\delta(M_{j'}, M_T)/2}\right)$$

and it follows that

$$\frac{P(M_j|\mathbf{y}, \mathbf{X})}{P(M_{j'}|\mathbf{y}, \mathbf{X})}\Big|M_T \to \begin{cases} 0, & \text{if } \delta(M_{j'}, M_T) < \delta(M_j, M_T), \\ \infty, & \text{if } \delta(M_{j'}, M_T) > \delta(M_j, M_T). \end{cases}$$

Thus, the model that is closer to $M_T$ is always preferred.

4. If the distance from both models to the true one is the same, that is, $\delta(M_{j'}, M_T) = \delta(M_j, M_T)$, then the limiting behavior of the quotient of posterior model probabilities only depends on the number of covariates of the models. We have that

$$\frac{P(M_j|\mathbf{y}, \mathbf{X})}{P(M_{j'}|\mathbf{y}, \mathbf{X})}\Big|M_T \to \begin{cases} 0, & \text{if } \delta(M_{j'}, M_T) = \delta(M_j, M_T) \text{ and } j' < j, \\ 1, & \text{if } \delta(M_{j'}, M_T) = \delta(M_j, M_T) \text{ and } j' = j, \\ \infty, & \text{if } \delta(M_{j'}, M_T) = \delta(M_j, M_T) \text{ and } j' > j. \end{cases}$$

(15)

When the true model is nested in $M_j$ and $M_{j'}$, so $\delta(M_{j'}, M_T) = \delta(M_j, M_T)$, (15) says that the smaller model is then preferred. Thus, the intrinsic Bayes procedure naturally leans toward a more parsimonious solution.

5. We also address the important point of what happens when the true model is a linear model but it does not belong to $\mathfrak{M}$. This happens when, for example, the true model includes some covariates or interactions among the existing or new ones not previously considered. From the preceding discussion it follows easily that the preference of the models in $\mathfrak{M}$ solely depends on their distances to the true model, regardless of whether the latter does or does not belong to the set of models we are considering.

Lastly, we note that implementation of the model selection procedure is best done with a stochastic search algorithm. As there are $2^{k-1}$ possible models, enumeration quickly becomes infeasible. We have implemented Metropolis-Hastings driven stochastic searches for both variable selection [Casella and Moreno (2006)] and changepoint problems [Girón, Moreno and Casella (2007)] with good results.

### APPENDIX: DERIVATION OF THE INTRINSIC BAYES FACTOR

Here we outline the calculations to justify the intrinsic Bayes factor of (7). For comparing the models in (6) with

$$\pi^I(\boldsymbol{\beta}_j, \sigma_j | \boldsymbol{\alpha}_i, \sigma_i) = \frac{2}{\pi \sigma_i (1 + \sigma_j^2 | \sigma_i^2)} N_j(\boldsymbol{\beta}_j | \tilde{\boldsymbol{\alpha}}_j, (\sigma_j^2 + \sigma_i^2) \mathbf{W}_j^{-1}),$$

$$\pi^I(\boldsymbol{\beta}_j, \sigma_j) = \int \pi^I(\beta_j, \sigma_j | \boldsymbol{\alpha}_i, \sigma_i) \pi^N(\boldsymbol{\alpha}_i, \sigma_i) \, d\boldsymbol{\alpha}_i \, d\sigma_i$$

and

$$\mathbf{W}_j^{-1} = \frac{n}{j+1} (\mathbf{X}_j' \mathbf{X}_j)^{-1},$$

the Bayes factor is given by (7).

The derivation of this expression is similar to that in Casella and Moreno (2006), but there different default priors were used and a generic $\mathbf{W}_j$ was derived. Here, we are using the reference prior $\pi^N(\boldsymbol{\eta}, \sigma) = c/\sigma$ instead, which seems to be a better choice as discussed in Girón et al. (2006), and thus we obtain a slightly different Bayes factor given by

$$BF_{ji}^n = \frac{2}{\pi} |\mathbf{X}_i' \mathbf{X}_i|^{1/2} (\mathbf{y}'(\mathbf{I}_n - \mathbf{H}_i)\mathbf{y})^{(n-i)/2} I_0,$$

where

$$I_0 = \int_0^{\pi/2} \frac{d\varphi}{|\mathbf{A}(\varphi)|^{1/2} |\mathbf{B}(\varphi)|^{1/2} E(\varphi)^{n-i}},$$

$$\mathbf{B}(\varphi) = \sin^2 \varphi \mathbf{I}_n + \mathbf{X}_j \mathbf{W}_j^{-1} \mathbf{X}_j',$$

$$\mathbf{A}(\varphi) = \mathbf{X}_i' \mathbf{B}(\varphi)^{-1} \mathbf{X}_i$$

and

$$E(\varphi) = \mathbf{y}'(\mathbf{B}(\varphi)^{-1} - \mathbf{B}(\varphi))^{-1} \mathbf{X}_i \mathbf{A}(\varphi)^{-1} \mathbf{X}_i' \mathbf{B}(\varphi)^{-1} \mathbf{y}.$$

Now, taking

$$\mathbf{W}_j^{-1} = \frac{n}{j+1} (\mathbf{X}_j' \mathbf{X}_j)^{-1}$$

we have, after some algebra, the following equalities:

(i)

$$\mathbf{B}(\varphi)^{-1} = \frac{1}{\sin^2 \varphi} \left( \mathbf{I}_n - \frac{n}{n + (j+1)\sin^2 \varphi} \mathbf{H}_j \right),$$

(ii)

$$\mathbf{A}(\varphi) = \frac{j+1}{n + (j+1)\sin^2 \varphi} \mathbf{X}_i' \mathbf{X}_i,$$

(iii)

$$\mathbf{X}_i \mathbf{A}(\varphi)^{-1} \mathbf{X}_i' = \frac{n + (j+1)\sin^2 \varphi}{j+1} \mathbf{H}_i,$$

(iv)

$$E(\varphi) = \frac{j+1}{n + (j+1)\sin^2 \varphi} \left( \frac{n}{(j+1)\sin^2 \varphi} RSS_j + RSS_i \right),$$

(v)

$$|\mathbf{A}(\varphi)| = \left( \frac{j+1}{n + (j+1)\sin^2 \varphi} \right)^i |\mathbf{X}_i' \mathbf{X}_i|,$$

(vi)

$$|\mathbf{B}(\varphi)| = (\sin^2 \varphi)^{n-j} \left( \frac{n + (j+1)\sin^2 \varphi}{j+1} \right)^j.$$

Plugging these values into $I_0$ and making some simplifications we get expression (7).

G. CASELLA
DEPARTMENT OF STATISTICS
UNIVERSITY OF FLORIDA
GAINESVILLE, FLORIDA 32611
USA
E-MAIL: casella@stat.ufl.edu

F. J. GIRÓN
M. LINA MARTÍNEZ
DEPARTMENT OF STATISTICS
UNIVERSITY OF MÁLAGA
MÁLAGA
SPAIN
E-MAIL: fj_giron@uma.es
        mlmartinez@uma.es

E. MORENO
DEPARTMENT OF STATISTICS
UNIVERSITY OF GRANADA
18071, GRANADA
SPAIN
E-MAIL: emoreno@ugr.es